\documentclass[12pt]{article}
\usepackage{amsfonts}
\usepackage{}
\usepackage{wasysym,amssymb,eufrak,indentfirst,graphicx,cite,amsthm,color}
\usepackage[bookmarksnumbered, colorlinks, plainpages]{hyperref}
\newtheorem{Theorem}{Theorem}[section]
\theoremstyle{definition}

\newtheorem{Lemma}[Theorem]{Lemma}

\newtheorem{Example}[Theorem]{Example}
\theoremstyle{definition}
\newtheorem{Remark}{\bf{Remark}}

\begin{document}
\baselineskip 17pt

\title{On finite groups factorized by $\sigma$-nilpotent subgroups\thanks{Research was supported by NSFC of China(No. 12001526) and Natural Science Foundation of Jiangsu Province, China (No. BK20200626).}}

\author{Zhenfeng Wu\\
{\small  School of Science, Jiangnan University}\\
{\small Wuxi, 214122, P.R. China}\\
{\small E-mail: zhfwu@jiangnan.edu.cn}\\ \\
Chi Zhang\thanks{Corresponding author}\\
{\small Department of Mathematics, China University of Mining and Technology}\\
{\small Xuzhou 221116, P.R. China}\\
{\small E-mail: zclqq32@cumt.edu.cn}}

\date{}
\maketitle

\begin{abstract}
Let $G$ be a finite group
and $\sigma=\{\sigma_{i}|i\in I\}$ be a partition of the set of all primes $\mathbb{P}$, that is, $\mathbb{P}=\bigcup_{i\in I}\sigma_{i}$ and $\sigma_{i}\cap \sigma_{j}=\emptyset$ for all $i\neq j$.
A chief factor $H/K$ of $G$ is said to be $\sigma$-central in $G$, if the semidirect product $(H/K)\rtimes(G/C_G(H/K))$ is a $\sigma_i$-group for some $i\in I$.
The group $G$ is said to be $\sigma$-nilpotent if either $G=1$ or every chief factor of $G$ is $\sigma$-central.
In this paper, we  study the properties of a finite group $G=AB$, factorized by two $\sigma$-nilpotent subgroups $A$ and $B$, and also generalize some known results.
\end{abstract}

\let\thefootnoteorig\thefootnote
\renewcommand{\thefootnote}{\empty}


 \let\thefootnote\thefootnoteorig

\section{Introduction}
Throughout this paper, all groups  are finite and $G$ always denotes a  group.
Moreover, $n$ is an integer, $\mathbb{P}$ is the set of all primes.
The symbol $\pi(n)$ denotes the set of all primes dividing $n$ and $\pi(G)=\pi(|G|)$,
the set of all primes dividing the order of $G$.
And $\pi$ always denotes a set of primes.
Following \cite{hall}, we say that a finite group $G$ possesses the following properties: $E_{\pi}$ if $G$ contains a Hall $\pi$-subgroup;
$C_{\pi}$ if $G$ enjoys $E_{\pi}$, and every two Hall $\pi$-subgroups of $G$ are conjugate;
$D_{\pi}$ if $G$ enjoys $C_{\pi}$, and every one of its $\pi$-subgroups is contained in some Hall $\pi$-subgroup of $G$.

In what follows, $\sigma=\{\sigma_{i}|i\in I\}$ is some partition of $\mathbb{P}$,
that is, $\mathbb{P}=\bigcup_{i\in I}\sigma_{i}$ and $\sigma_{i}\cap \sigma_{j}=\emptyset$ for all $i\neq j$.
$\Pi$ is always supposed to be a non-empty subset of the set $\sigma$ and $\Pi^{'}=\sigma\backslash \Pi$.
We write $\sigma(n)=\{\sigma_{i}|\sigma_{i}\cap \pi(n)\neq\emptyset\}$
and $\sigma(G)=\sigma(|G|)$.

Following \cite{AN3} and \cite{A},
$G$ is said to be $\sigma$-primary if $G=1$ or $|\sigma(G)|=1$.
A chief factor $H/K$ of $G$ is said to be $\sigma$-central in $G$ if the semidirect product $H/K\rtimes G/C_G(H/K)$ is $\sigma$-primary.
A set $\mathcal {H}$ of subgroups of $G$ is said to be a complete Hall $\sigma$-set of $G$ if every non-identity member of $\mathcal {H}$ is a Hall $\sigma_{i}$-subgroup of $G$ for some $i$ and $\mathcal {H}$ contains exactly one Hall $\sigma_{i}$-subgroup of $G$ for every $\sigma_{i}\in \sigma(G)$.

Recall that $G$ is called \cite{AN3}: (i) $\sigma$-nilpotent if either $G=1$ or every chief factor of $G$ is $\sigma$-central in $G$,
(ii) $\sigma$-soluble if either $G=1$ or every chief factor of $G$ is $\sigma$-primary.
We use $\mathfrak{N}_\sigma$ and $\mathfrak{S}_\sigma$ to denote the classes of all $\sigma$-nilpotent and $\sigma$-soluble groups, respectively.
Following \cite{AN3} we call the product of all normal $\sigma$-nilpotent subgroups of $G$ the $\sigma$-Fitting subgroup of $G$ and denote it by $F_\sigma(G)$.
Clearly, $F_\sigma(G)$ is also $\sigma$-nilpotent (see the below Lemma \ref{normalnil}).

It has been established earlier by the work of Wielandt \cite{wie} and Kegel \cite{kegel} that the product of two finite nilpotent groups is soluble (see for example \cite[Theorem 2.4.3]{bf}).
Moreover, the structure of finite groups factorized by two nilpotent subgroups has been investigated by several authors (see \cite[Chapter 2]{bf}), and properties of such groups have also been discovered by many authors (see for example \cite{pa,hh1,hh2,pg,se,js}).
On the other hand, very little is known about the properties of a product of two $\sigma$-nilpotent subgroups.
So the aim of this paper is to extend the knowledge of properties of such products by some $\sigma$-nilpotent subgroups.

Robinson and Stonehewer \cite[Theorem 2]{rs} have shown that if the group $G=AB$ is the product of two abelian subgroups $A$ and $B$, then every chief factor of $G$ either is centralized by $A$ or $B$.
Moreover, Stonehewer \cite[Theorem 1]{se} have proved that if $G=AB$ is the product of two nilpotent subgroups $A$ and $B$, then for every minimal normal subgroup $N$ of $G$ one of the subgroups $AN$ and $BN$ is nilpotent.
Our first main theorem generalizes these results to a finite group factorized by two $\sigma$-nilpotent subgroups.

\begin{Theorem}\label{T1}
Let $G$ be the product of two $\sigma$-nilpotent subgroups $A$ and $B$ and let $N$ be a minimal normal subgroup of $G$.
Assume also that $G$ is $\sigma$-soluble.
Then $AN$ or $BN$ is $\sigma$-nilpotent.
\end{Theorem}

Note that, the product of two finite nilpotent groups is always soluble by the well-known theorem of Wielandt \cite{wie} and Kegel \cite{kegel}.
However, the product of two finite $\sigma$-nilpotent groups is not necessarily $\sigma$-soluble.
For example, let $G=A_5$ be the alternating group with degree 5, and let $\sigma_1=\{2,3\}$, $\sigma_2=\{5\}$ and $\sigma_3=\{2,3,5\}'$.
Then it is clear that $A_4$ is $\sigma$-nilpotent and the Sylow $5$-subgroup $P$ of $G$ is also $\sigma$-nilpotent.
Then $G=A_4P$ is the product of two $\sigma$-nilpotent groups, but clear $G$ is not $\sigma$-soluble.

It was also proved by Robinson and Stonehewer \cite[Theorem 1]{rs} that if a group $G$ has three abelian subgroups $A,B$ and $C$ such that $G=AB=BC=CA$, then every chief factor of $G$ is central in $G$.
This result has been extended to nilpotent subgroups, that is, Kegel \cite{kegel2} have shown that if the finite group $G=AB=BC=CA$ is the product of three nilpotent subgroups $A,B$ and $C$, then $G$ is nilpotent.
So we have the following theorem.

\begin{Theorem}\label{T2}
Let $G$ be a finite group with $\sigma$-nilpotent subgroups $A,B$ and $C$ such that $G=AB=BC=CA$.
Assume also that $G$ satisfies $D_{\sigma_i}$ for some $\sigma_i\in \sigma(G)$.
Then $G$ is $\sigma$-nilpotent.
\end{Theorem}

Later Cossey and Stonehewer studied the structure of the product of two nilpotent subgroups in which its Fitting subgroup is a $p$-group for some prime $p$, that is, Cossey and Stonehewer \cite[Theorem 1]{js} have proved that: Let $G$ be a soluble group for which $F(G)$ is a $p$-group (for some prime $p$).
Then $G$ is the product of two nilpotent subgroups if and only if it has a nilpotent Hall $p'$-subgroup.
Base on this fact, we can study the case of $\sigma$-nilpotent subgroups.
So we have the following theorem.

\begin{Theorem}\label{T3}
Let $G$ be a $\sigma$-soluble group with $F_\sigma(G)$ is a $\sigma_i$-group for some $\sigma_i\in\sigma(G)$.
Then $G=AB$, where $A,B$ are $\sigma$-nilpotent subgroups of $G$ if and only if $G$ has a $\sigma$-nilpotent Hall $\sigma_i'$-subgroup.
\end{Theorem}

\vskip 0.25cm
All unexplained terminologies and notations are standard,
as in \cite{KT} and \cite{AN3}.

\section{Preliminaries}
\begin{Lemma}\label{soluble}\textup(See {\cite[Lemma 2.1]{AN2}}) The class $\mathfrak{S}_\sigma$ is closed under taking direct products, homomorphic images and subgroups.
Moreover, any extension of the $\sigma$-soluble group by a $\sigma$-soluble group is a $\sigma$-soluble group as well.
\end{Lemma}

Recall that the product of all normal $\sigma$-soluble subgroups of $G$ called the $\sigma$-radical \cite{AN2} of $G$ and denote it by $R_\sigma(G)$.

\begin{Lemma}\label{D}\textup(See {\cite[Theorem B ]{AN2}}) Let $R=R_\sigma(G)$ be the $\sigma$-radical of $G$.
Then $G$ is $\sigma$-soluble if and only if for any $\Pi$ the following hold: $G$ has a Hall $\Pi$-subgroup $E$, every $\Pi$-subgroup of $G$ is contained in some conjugate of $E$ and $E$ $R$-permutes with every Sylow subgroup of $G$.
\end{Lemma}

\begin{Lemma}\label{product}\textup(See {\cite[Lemma ]{pa}}) Let $G=AB$ be a finite group satisfying $D_{\pi}$.
If the subgroups $A$ and $B$ of $G$ possess normal Hall $\pi$-subgroups $A_{\pi}$ and $B_{\pi}$, respectively,
then $A_{\pi}B_{\pi}=B_{\pi}A_{\pi}$ is a Hall $\pi$-subgroup of $G$ and $[A_{\pi},B_{\pi}]\leq O_{\pi}(G)$.
In particular, if $O_{\pi}(G)=1$, then $[A_{\pi}^G, B_{\pi}^G]=1$.
\end{Lemma}

\begin{Lemma}\label{proper}\textup(See {\cite[Satz 3]{kegel}}) Let $A$ and $B$ be the subgroups of $G$ such that $G\neq AB$ and $AB^x=B^xA$ for all $x\in G$.
Then $G$ has a proper normal subgroup $N$ such that either $A\leq N$ or $B\leq N$.
\end{Lemma}

\begin{Lemma}\label{contain}\textup(See {\cite[Chapter 2, Lemma 4.4]{Guo3}}) Let $G=AB$ be the product of the subgroups $A$ and $B$.
If $L$ is a normal subgroup of $A$ and $L\leq B$, then $L\leq B_G$.
\end{Lemma}

\begin{Lemma}\label{Dproperty}\textup(See {\cite[Lemma 2.1]{revin}}) Let $G$ be a finite group and $A$ its normal subgroup.
If $G$ satisfies $D_{\pi}$, then $G/A$ satisfies $D_{\pi}$.
\end{Lemma}

\begin{Lemma}\label{nilpotent}\textup(See {\cite[Proposition 3.4]{scie}})
Any two of the following conditions are equivalent:

$(1)$ $G$ is $\sigma$-nilpotent.

$(2)$ $G=H_1\times\cdots\times H_t$, where $\{H_1,...,H_t\}$ is a complete Hall $\sigma$-set of $G$.

$(3)$ Every subgroup of $G$ is $\sigma$-subnormal in $G$.
\end{Lemma}

\begin{Lemma}\label{normalnil}\textup(See {\cite[Corollary 2.4]{AN3}})
If $A$ and $B$ are normal $\sigma$-nilpotent subgroups of $G$, then $AB$ is $\sigma$-nilpotent.
\end{Lemma}

\begin{Lemma}\label{con}\textup(See {\cite[Chapter 1, Lemma 4.14]{Guo3}})
Let $A$ and $B$ be proper subgroups of $G$ such that $G=AB$.
Then:

(1) $G_p=A_pB_p$ for some Sylow $p$-subgroups $G_p$, $A_p$ and $B_p$ of $G$, $A, B$, respectively;

(2) $G=A^xB$ and $G\neq AA^x$ for all $x\in G$.

\end{Lemma}

\section{Proofs of Theorem \ref{T1}, \ref{T2} and \ref{T3}}

{\bf Proof of Theorem \ref{T1}}. Assume that the result is false and let $G$ be a counterexample of minimal order.
Since $G$ is $\sigma$-soluble and $N$ is a minimal normal subgroup of $G$, we have that $N<G$ and $N$ is a $\sigma_i$-group for some $i$.
Now without loss of generality we may assume that $N$ is a $\sigma_1$-group, and let $\Pi=\sigma_1'$.
Then we proceed the proof by the following steps.

$(1)$ $N$ is the unique minimal normal subgroup of $G$.

If not, then let $L$ be a minimal normal subgroup of $G$ different from $N$.
Since $G/L=(AL/L)(BL/L)$ and $AL/L\cong A/(A\cap L)$, $BL/L\cong B/(B\cap L)$ are all $\sigma$-nilpotent,
we have that $G/L$ satisfies the hypothesis.
Hence by the choice of $G$, we see that $(AL/L)(NL/L)$ or $(BL/L)(NL/L)$ is $\sigma$-nilpotent.
If $(AL/L)(NL/L)$ is $\sigma$-nilpotent, then
\begin{center}
 $AN/(AN\cap L)\cong ANL/L=(AL/L)(NL/L)$
\end{center} is $\sigma$-nilpotent.
Clearly, $AN/N\cong A/(A\cap N)$ is $\sigma$-nilpotent.
So
\begin{center}
 $AN=(AN)/(AN\cap L\cap N)\lesssim (AN/N)\times (AN/AN\cap L)$
\end{center}
is $\sigma$-nilpotent, a contradiction.
Similarly, if $(BL/L)(NL/L)$ is $\sigma$-nilpotent, then we obtain that $BN$ is also $\sigma$-nilpotent, a contradiction.
Hence $(1)$ holds.

$(2)$ $O_{\Pi}(G)=1$.

This is clear from $(1)$.

Since $A$ and $B$ are $\sigma$-nilpotent, we get that $A$ and $B$ have normal Hall $\Pi$-subgroup $A_{\Pi}$ and $B_{\Pi}$, respectively.

$(3)$ $A_{\Pi}\neq1$ and $B_{\Pi}\neq1$.

If $A_{\Pi}=1$, then $A$ is a $\sigma_1$-group.
This implies that $AN$ is also a $\sigma_1$-group, and so is $\sigma$-nilpotent, a contradiction.
Hence $A_{\Pi}\neq1$.
Similarly, $B_{\Pi}\neq1$.
Hence $(3)$ holds.

$(4)$ Final contradiction.

Since $G$ is $\sigma$-soluble, we have that $G$ satisfies $D_{\Pi}$ be Lemma \ref{D}.
Clearly, $A_{\Pi}\unlhd A$ and $B_{\Pi}\unlhd B$.
Then by $(2)$ and Lemma \ref{product}, we see that $[A_{\Pi}^G,B_{\Pi}^G]=1$, and so $A_{\Pi}^G\leq C_G(B_{\Pi}^G)$ and $B_{\Pi}^G\leq C_G(A_{\Pi}^G)$.
From $(1)$ and $(3)$ we get that $N\leq A_{\Pi}^G\cap B_{\Pi}^G$.
Hence $A_{\Pi}^G\leq C_G(N)$ and $B_{\Pi}^G\leq C_G(N)$.
It is obvious that $A_{\Pi}$ is also a Hall $\Pi$-subgroup of $AN$.
Since $A$ is $\sigma$-nilpotent, and $A_{\Pi}\leq A_{\Pi}^G\leq C_G(N)$, we have that $A_{\Pi}$ is the normal Hall $\Pi$-subgroup of $AN$.
Clearly, $A_{\sigma_1}N$ is a Hall $\sigma_1$-subgroup of $AN$ and $A_{\sigma_1}N$ is normal in $AN$.
Then $AN=A_{\Pi}\times A_{\sigma_1}N$ and so $AN$ is $\sigma$-nilpotent, a contradiction.
Similarly, $BN$ is also $\sigma$-nilpotent, a contradiction.
The final contradiction completes the theorem.

\vskip 0.25cm
{\bf Proof of Theorem \ref{T2}}. First of all, we show that $G$ is $\sigma$-soluble.
Assume that the assertion is false and let $G$ be a counterexample of minimal order.
Without loss of generality we may assume that $\sigma_i\in \sigma(B)$.
Then by Lemma \ref{nilpotent}, we see that $B$ has a nonidentity normal Hall $\sigma_i$-subgroup, denote it by $B_i$.
If $\sigma_i\notin\sigma(A)$, then $|G:C|=|A:A\cap C|$ is a $\sigma_i'$-number for $G=AC$.
Since $G=BC$, we see that $|B:B\cap C|=|BC:C|=|G:C|$ is a $\sigma_i'$-number.
Hence $B_i\leq B\cap C$ for $B$ is $\sigma$-nilpotent.
Therefore by Lemma \ref{contain}, we have that $1\neq B_i\leq C_G$.
Now let $R$ be a minimal normal subgroup of $G$ which is contained in $C_G$.
Then $R$ is $\sigma$-soluble.
Clearly, by Lemma \ref{Dproperty} we see that $G/R$ satisfies $D_{\sigma_i}$.
Since
\begin{center}
 $G/R=(AR/R)(BR/R)=(BR/R)(CR/R)=(CR/R)(AR/R)$
\end{center}
and
\begin{center}
 $AR/R\cong A/A\cap R$, $BR/R\cong B/B\cap R$ and $CR/R\cong C/C\cap R$
\end{center}
are $\sigma$-nilpotent, we have that $G/R$ satisfies the hypothesis.
Hence by the choice of $G$, we obtain that $G/R$ is $\sigma$-soluble.
Then by Lemma \ref{soluble} shows that $G$ is $\sigma$-soluble for $R$ is $\sigma$-soluble, a contradiction.
This contradiction shows that $\sigma_i\in\sigma(A)$.
Since $A$ is $\sigma$-nilpotent, we have that $A$ has a nonidentity Hall $\sigma_i$-subgroup by Lemma \ref{nilpotent}, and denote it by $A_i$.
Then by Lemma \ref{product}, we see that $A_iB_i=B_iA_i$ is a Hall $\sigma_i$-subgroup of $G$.
It is obvious that $A_iB_i<G$ and for any element $x\in G$, by Lemma \ref{con} $G=AB^x$.
And clearly, $B^x$ is also $\sigma$-nilpotent and $B_i^x$ is the normal Hall $\sigma_i$-subgroup of $B^x$.
Hence by Lemma \ref{product} again, we have that $A_iB_i^x=B_i^xA_i$ is a Hall $\sigma_i$-subgroup of $G$.
So by Lemma \ref{proper}, we obtain that there exists a proper normal subgroup $N$ of $G$ such that $A_i\leq N$ or $B_i\leq N$.
This means that $G$ is not a simple group.
Let $R$ be a minimal normal subgroup of $G$.
Then $G/R$ satisfies the hypothesis by Lemma \ref{Dproperty}.
Hence by the choice of $G$ we have that $G/R$ is $\sigma$-soluble.
Therefore $R$ is the unique minimal normal subgroup of $G$ and $O_{\sigma_i}(G)=1$.
So by Lemma \ref{product}, we see that $[A_{\sigma_i}^G,B_{\sigma_i}^G]=1$, that is, $A_{\sigma_i}^G\leq C_G(B_{\sigma_i}^G)$.
Since $A_{\sigma_i}\neq 1$ and $B_{\sigma_i}\neq 1$ and $R$ is the unique minimal normal subgroup of $G$, we obtain that $R\leq A_{\sigma_i}^G\cap B_{\sigma_i}^G$.
So $R\leq A_{\sigma_i}^G\leq C_G(B_{\sigma_i}^G)\leq C_G(R)$, which shows that $R$ is abelian.
Hence we also get $G$ is $\sigma$-soluble for $G/R$ is $\sigma$-soluble, a contradiction.
This contradiction shows that $G$ is $\sigma$-soluble.

Now let $N$ be a minimal normal subgroup of $G$.
Then since $G$ is $\sigma$-soluble, we have that $N$ is a $\sigma_j$-group for some $\sigma_j\in\sigma(G)$.
If $N=1$, then $G$ is a simple $\sigma$-soluble, and so $G$ is $\sigma_i$-group, that is, $G$ is $\sigma$-nilpotent.
So we can always assume that $N\neq1$.
Now without loss of generality we may assume that $N$ is a $\sigma_1$-group.
Then we can proceed to prove that $G$ is $\sigma$-nilpotent by using induction on $|G|$.
Since
\begin{center}
   $G/N=(AN/N)(BN/N)=(BN/N)(CN/N)=(CN/N)(AN/N)$
\end{center}
and
\begin{center}
 $AN/N\cong A/A\cap N$, $BN/N\cong B/B\cap N$ and $CN/N\cong C/C\cap N$
\end{center}
are $\sigma$-nilpotent, and by Lemma \ref{Dproperty} $G/N$ satisfies $D_{\sigma_i}$, we have that $G/N$ satisfies the hypothesis.
Hence by induction, we have that $G/N$ is $\sigma$-nilpotent.
If $G$ has another minimal normal subgroup $R$ of $G$ which is different from $N$, then by the same discussion as above, we also get $G/R$ is $\sigma$-nilpotent.
Hence $G=G/(N\cap R)\lesssim G/N\times G/R$ is $\sigma$-nilpotent.

Now we can only consider that $G$ has a unique minimal normal subgroup $N$ of $G$.
By the discuss as above, we see that $G/N$ is $\sigma$-nilpotent.
So if $N\leq \Phi(G)$, then $G$ is also $\sigma$-nilpotent.
Hence we can assume that $N\nleq\Phi(G)$.
Since by Theorem \ref{T1} and $G$ is $\sigma$-soluble, we have that at least two of $AN$, $BN$ and $CN$ are $\sigma$-nilpotent.
So we can assume that $AN$ and $BN$ are $\sigma$-nilpotent.
This implies that $A_{\sigma_1'}\leq C_G(N)$ and $B_{\sigma_1'}\leq C_G(N)$ by Lemma \ref{nilpotent}.
If $A_{\sigma_1'}\neq1$, then $A_{\sigma_1'}\leq C_G(N)$, it follows that $C_G(N)\neq1$.
But $N$ is the unique minimal normal subgroup of $G$, we have that $N\leq C_G(N)$, that is, $N$ is an elementary abelian $p$-group.
Since $N\nleq\Phi(G)$, there exists a maximal subgroup $M$ of $G$ such that $G=NM=C_G(N)M$.
Clearly, $N\cap M\unlhd G$ and $C_G(N)\cap M\unlhd G$.
This shows that $N\cap M=1$ and $C_G(N)\cap M=1$ for $N$ is the unique minimal normal subgroup of $G$.
So $N=C_G(N)$.
But $A_{\sigma_1'}\leq C_G(N)=N$ and $N$ is a $\sigma_1$-group, we have that $A_{\sigma_1'}=1$, a contradiction.
This contradiction shows that $A_{\sigma_1'}=1$.
If we assume that $B_{\sigma_1'}\neq1$, then we can get a contradiction like above.
So this contradiction shows that $B_{\sigma_1'}=1$.
Hence $A_{\sigma_1'}=1$ and $B_{\sigma_1'}=1$.
But $A$ and $B$ are $\sigma$-nilpotent, we have that $A$ and $B$ are all $\sigma_1$-group, which means that $G=AB$ is also a $\sigma_1$-group.
So $G$ is $\sigma$-nilpotent.
Thus we have proved that $G$ is $\sigma$-nilpotent.

\vskip 0.25cm
{\bf Proof of Theorem \ref{T3}}.
On one hand, if $G$ has a $\sigma$-nilpotent Hall $\sigma_i'$-subgroup $K$.
Then since $G$ is $\sigma$-soluble, we have that there exists a Hall $\sigma_i$-subgroup $H$ of $G$.
Clearly, $G=HK$ and $H,K$ are all $\sigma$-nilpotent.

On the other hand, suppose that $G=AB$ with $A,B$ are $\sigma$-nilpotent subgroups of $G$ and $F_\sigma(G)$ is a $\sigma_i$-group.
First, we claim that $O_{\sigma_i'}(G)=1$.
In fact, if $O_{\sigma_i'}(G)\neq1$, then let $N$ be a minimal normal subgroup of $G$ contained in $O_{\sigma_i'}(G)$.
So $N$ is a $\sigma_j$-group for some $\sigma_{j}\in\sigma(G)\backslash\sigma_i$, and so $N\leq F_\sigma(G)$.
This contradicts to the fact that $F_\sigma(G)$ is a $\sigma_i$-group.
Hence we have $O_{\sigma_i'}(G)=1$.
Since $A$ and $B$ are $\sigma$-nilpotent and by Lemma \ref{nilpotent}, we have that $A$ and $B$ have normal Hall $\sigma_i'$-subgroup $A_{\sigma_i'}$ and $B_{\sigma_i'}$, respectively.
By Lemma \ref{D}, we see that $G$ satisfies $D_{\sigma_i'}$ for $G$ is $\sigma$-soluble.
Hence by Lemma \ref{product}, we obtain that $G_{\sigma_i'}= A_{\sigma_i'}B_{\sigma_i'}$ is the Hall $\sigma_i'$-subgroup of $G$ and $[A_{\sigma_i'}^G,B_{\sigma_i'}^G]=1$.
If $A_{\sigma_i'}=1$ or $B_{\sigma_i'}=1$, then $G_{\sigma_i'}=B_{\sigma_i'}\leq B$ or $G_{\sigma_i'}=A_{\sigma_i'}\leq A$, which all shows that $G_{\sigma_i'}$ is $\sigma$-nilpotent for $A$ and $B$ are $\sigma$-nilpotent.
If $A_{\sigma_i'}\neq1$ and $B_{\sigma_i'}\neq1$, then since $[A_{\sigma_i'}^G,B_{\sigma_i'}^G]=1$, we have that $A_{\sigma_i'}\leq A_{\sigma_i'}^G\leq C_G(B_{\sigma_i'}^G)\leq C_G(B_{\sigma_i'})$.
This implies that $A_{\sigma_i'}\unlhd G_{\sigma_i'}$ and $B_{\sigma_i'}\unlhd G_{\sigma_i'}$.
Thus by Lemma \ref{normalnil}, we have that $G_{\sigma_i'}$ is $\sigma$-nilpotent.
Therefore in any case, we always get that $G$ has a $\sigma$-nilpotent Hall $\sigma_i'$-subgroup.

\begin{Remark}
Note that the  Theorem \ref{T1} does not hold if "$G$ is $\sigma$-soluble" is replaced by "$G$ satisfies $D_{\sigma_i}$, for every $\sigma_i\in\sigma(G)$".
We will give a counterexample to show this.
\end{Remark}

\begin{Example}
Let $G=PSL_2(7)$ be the projective special linear group of degree 2 over $F_7$, a simple group with order $168=2^3\cdot3\cdot7$.
Let $\pi=\{3,7\}$.
Then by \cite[Condition III]{vdovin}, we see that $G$ satisfies $D_{\pi}$.
Now let $H$ be a Hall $\pi$-subgroup of $G$, and $P$ be a Sylow $2$-subgroup of $G$, and
$\sigma_1=\{2\}$, $\sigma_2=\{3,7\}$ and $\sigma_3=\{2,3,7\}'$.
So $P$ is a Hall $\sigma_1$-subgroup of $G$ and $H$ is a Hall $\sigma_2$-subgroup of $G$ and $G=HP$.
Clearly, $G$ satisfies $D_{\sigma_1}$ and $D_{\sigma_2}$.
The  Theorem \ref{T1} does not hold for $G$.
However $G$ is not $\sigma$-soluble.
\end{Example}

\end{document}